\documentclass[12pt]{article}

\usepackage[french]{babel}
\usepackage[utf8]{inputenc}  
\usepackage[T1]{fontenc}
\usepackage{amsfonts,a4wide,amsmath,amssymb,amsthm}
\usepackage[all]{xy}    
\usepackage{mathtools}  
\usepackage{multirow}    
\usepackage{csquotes}
 \usepackage{relsize}
\usepackage{tikz}   
\usepackage{hyperref}      
\hypersetup{
     colorlinks=true, 
     breaklinks=true, 
     urlcolor= blue,  
     linkcolor= black, 
     linktoc	=all,	
     citecolor=red,	
     filecolor=black,	
     bookmarksopen=true,            
     pdftitle={La méthode de Runge pour des variétés modulaires}, 
     pdfauthor={Samuel Le Fourn},     
}

\allowdisplaybreaks

\theoremstyle{plain}
\newtheorem{thm}{Théorème}

\newtheorem*{thmsansnom}{Théorème}

\theoremstyle{definition}

\theoremstyle{remark}
 
\newcommand{\gP}{{\mathfrak{P}}}

\newcommand{\Acal}{{\mathcal A}}

\newcommand{\Dcal}{{\mathcal D}}
\newcommand{\Ecal}{{\mathcal E}}
\newcommand{\Fcal}{{\mathcal F}}

\newcommand{\Hcal}{{\mathcal H}}

\newcommand{\Ocal}{{\mathcal O}}

\newcommand{\Ucal}{{\mathcal U}}

\newcommand{\Xcal}{{\mathcal X}}

\newcommand{\Z}{{\mathbb{Z}}}
\newcommand{\Q}{{\mathbb{Q}}}

\newcommand{\C}{{\mathbb{C}}}
\renewcommand{\P}{\mathbb{P}}

\newcommand{\Gal}{\operatorname{Gal}}

\renewcommand{\Im}{\operatorname{Im}}

\newcommand{\fonction}[5]{\begin{array}{c|ccl}           
#1: & #2 & \longrightarrow & #3 \\
    & #4 & \longmapsto & #5 \end{array}}

\newcommand{\quot}[2]                                    
{\raisebox{.6ex}{\newline$#1$}\!/\!\raisebox{-.6ex}{$#2$}}

\newcommand{\Kb}{\overline{K}}
\newcommand{\Qb}{\overline{\Q}}

\title{Sur la méthode de Runge et les points entiers de certaines variétés modulaires de Siegel}
\author{Samuel Le Fourn\footnote{Email : \url {samuel.le_fourn@ens-lyon.fr}  }
\\ ENS de Lyon}

\date\today

\begin{document}
\maketitle

\begin{abstract}
	Dans cette note, nous annonçons des résultats sur les points entiers de certaines variétés modulaires basés sur une généralisation de la méthode dite de Runge en dimension supérieure qui sera expliquée dans un premier temps. En particulier, nous obtenons un résultat explicite dans le cas de la variété modulaire de Siegel $A_2(2)$.
\end{abstract}
\section*{Introduction}

Dans cette note, nous décrivons les principaux résultats de \cite{LeFourn3} qui visent à établir de nouveaux éclairages sur le comportement des points entiers de variétés modulaires en dimension supérieure. 


Les questions de points entiers et rationnels de courbes modulaires jouent naturellement un grand rôle dans les recherches actuelles sur les propriétés des courbes elliptiques car ce sont leurs espaces de modules. Par exemple, Bilu et Parent \cite{BiluParent11} ont résolu le cas dit \og de Cartan déployé \fg{} du problème d'uniformité de Serre pour le corps $\Q$ en prouvant que  pour $p$ premier assez grand, les courbes modulaires notées $X_{\rm{split}}(p)$ n'ont pas de points rationnels non triviaux. Cette note est exclusivement consacrée à la \textit{méthode de Runge} : c'est un des trois outils de la preuve de Bilu et Parent et nous allons l'introduire plus bas. 

Avant de présenter la méthode, donnons-en directement une application effective (très simplifiée) obtenue par nos résultats. 
\begin{thmsansnom}
\hspace*{\fill}

	Soit $A$ une surface abélienne principalement polarisée sur un corps de nombres $K$ (rationnel ou quadratique imaginaire).
	Supposons que toute la 2-torsion de $A$ est définie sur $K$, et que $A$ a potentiellement bonne réduction en toute place finie de $K$.
	
	Alors, si la réduction semistable de $A$ est un produit de courbes elliptiques en au plus 3 places de $K$, la hauteur de Faltings stable de $A$ admet la borne 
	\[
	h_\Fcal(A) \leq 1070.
	\]
\end{thmsansnom}

Ainsi, le but de \cite{LeFourn3} est d'adapter la méthode de Runge en dimension supérieure, en particulier pour traiter non plus de courbes modulaires mais de variétés modulaires (par exemple celles de Siegel). Cette note s'articule autour des trois résultats principaux de l'article : le premier donne une généralisation de la méthode de Runge en dimension supérieure (partie \ref{secRungemethod}), conçue pour une certaine souplesse d'utilisation. Le second en déduit un résultat de finitude de points entiers \og à la Runge \fg{} sur les variétés modulaires de Siegel $A_2(n)$ (partie \ref{secvarmod}). Enfin, le troisième réalise explicitement la méthode dans le cas $n=2$, en tant que preuve de principe de son effectivité (partie \ref{seccasdeux}), avec notamment pour cas particulier le théorème ci-dessus. 

\section{La méthode de Runge}
\label{secRungemethod}
La méthode de Runge est la réalisation pratique de ce théorème (dû à Bombieri et inspiré par un résultat plus vieux de Runge)

\begin{thmsansnom}[Runge-Bombieri, Théorème 9.6.6 de \cite{BombieriGubler}]
\hspace*{\fill}

Soit $C$ une courbe projective lisse sur un corps de nombres $K$ et $\phi \in K(C)$ non constante. Pour une extension $L$ de $K$, $S_L$ peut désigner tout ensemble fini de places de $L$ contenant les places archimédiennes et on note $r_L$ le nombre d'orbites des pôles de $\phi$ par $\Gal(\Kb/L)$. Un couple $(L,S_L)$ vérifie la \og condition de Runge \fg{} si 
\begin{equation}
\label{eqcondRunge}
|S_L|<r_L.
\end{equation}
Alors, les points $P$ de $C(L)$ tels que $\phi(P)$ est $S_L$-entier (i.e. $\phi(P) \in \Ocal_{L,S_L}$) sont en nombre fini. En fait, on a même le résultat \og uniforme en les couples \fg{}
\begin{equation}
\label{eqRungefini}
	\bigcup_{\substack{(L,S_L) \\ |S_L|<r_L}} \{ P \in C(L) \, | \, \phi(P) \in \Ocal_{L,S_L} \} \textrm{  est fini.} 
\end{equation}
\end{thmsansnom}

Parmi les nombreuses méthodes de recherche de points entiers sur les courbes, celle-ci a deux avantages notables. Le premier est l'uniformité du résultat en les paires réalisant la condition de Runge : la plupart des méthodes demandent le choix d'une paire $(L,S_L)$ quelconque et donnent une borne sur la hauteur qui dépend de ce choix. Le second est l'effectivité de la méthode : si on connaît suffisamment bien les fonctions auxiliaires intervenant dans le théorème, la finitude s'écrit comme une borne effective sur la hauteur des points $\phi(P)$.

Le principe de la preuve est le suivant. Soit un couple $(L,S_L)$ fixé et $P \in C(L)$ tel que $\phi(P) \in \Ocal_{L,S_L}$. Le point $P$ est $v$-adiquement loin de tous les pôles pour toute place $v$ de $L$ n'appartenant pas à $S_L$, comme $|\phi(P)|_v \leq 1$. D'autre part, pour chaque $v \in S_L$, il ne peut s'approcher près que d'une orbite de pôles car celles-ci sont disjointes. Sous la condition de Runge, il reste donc une orbite $O$ de pôles qui est $v$-adiquement loin de $P$ quel que soit notre place $v$ sur $L$. En choisissant une fonction auxiliaire $g_O \in L(C)$ dont les pôles sont exactement $O$, la hauteur de Weil $h$ de $g_O(P)$ est alors bornée d'où la finitude de l'ensemble. En pratique, on peut même directement établir une borne absolue sur la hauteur de $\phi(P)$ et cela implique l'uniformité en les couples $(L,S_L)$.

Pour les courbes modulaires en général, Bilu et Parent ont exécuté cette méthode dans \cite{BiluParent09} avec la fonction $j$-invariant pour obtenir des bornes explicites sur les hauteurs des points entiers lorsque la condition de Runge est vérifiée. Par exemple, ils ont obtenu que si $P$ est un point de $X_{\rm{split}}(p)(\Q)$, alors $h(j(P)) = O (\sqrt{p})$.

Avant de discuter la généralisation, expliquons les données équivalentes à $(C,\phi)$ en dimension supérieure sous des hypothèses simplifiées. $C$ sera remplacé par un schéma projectif lisse $X$ sur $K$ et les pôles de $\phi$ correspondent à des diviseurs effectifs $D_1, \cdots, D_r$ sur $X$ (on suppose pour simplifier qu'ils sont définis sur $K$ et pas sur une extension), d'union $D$. L'intégralité est plus facilement formulée à partir d'un modèle projectif lisse $\Xcal$ sur $\Ocal_K$, avec $\Dcal$ la fermeture de Zariski de $D$ dans $\Xcal$, il s'agit alors de prouver la finitude de $(\Xcal \backslash \Dcal) (\Ocal_{L,S_L})$ sous une condition de type \eqref{eqcondRunge}. Soit donc $P \in (\Xcal \backslash \Dcal) (\Ocal_{L,S_L})$. Pour les places $v \notin S_L$, le point $P$ est bien $v$-adiquement loin de chaque $D_i$. Par contre, pour $v \in S_L$, il peut être proche de plusieurs $D_i$ en même temps car ils ne sont pas forcément disjoints... En fait, en notant $m$ le nombre maximal tels que l'intersection de $m$ diviseurs distincts parmi $D_1, \cdots, D_r$ est non vide, $P$ peut être proche de $m$ diviseurs $D_i$ en même temps, ce qui donne la condition de Runge multidimensionnelle
\begin{equation}
 \label{eqcondRungemulti}
	m |S_L|<r.
\end{equation}
Par ailleurs, la finitude du nombre de points $v$-adiquement loin pour toute place $v$ d'un certain diviseur revient à la propriété de Northcott sur la hauteur associée à ce diviseur, et on a donc besoin que chaque $D_i$ soit ample. Sous cette condition, l'argument ci-dessus (qui dissimule bien des complications techniques) montre que 
\begin{equation}
\label{eqRungemulti}
\bigcup_{\substack{(L,S_L) \\ m |S_L|<r}} (\Xcal \backslash \Dcal) (\Ocal_{L,S_L}) \textrm{ est fini.}
\end{equation}
Ce théorème, obtenu par Levin \cite{Levin08}, est tout à fait adapté pour certaines variétés, mais il ne semble malheureusement pas applicable aux variétés modulaires qui nous intéressent. La raison principale (mis à part l'hypothèse d'amplitude qui paraît inévitable) est que le nombre $m$ n'est pas assez petit par rapport à $r$ pour que la condition  \eqref{eqcondRungemulti} soit satisfaisable (n'oublions pas que $S_L$ contient toutes les places archimédiennes et plus généralement toutes les places de mauvaise définition du modèle $\Xcal$). Nous avons donc conçu une version plus souple de ce résultat permettant son application à d'autres variétés, basée sur la notion de \og voisinage tubulaire \fg{}. 

L'idée est qu'on suppose en plus que le point $P$ n'est pas trop près des points où les diviseurs $D_1, \cdots, D_r$ s'intersectent. Pour cela, on choisit un fermé $Y$ de $X$ et on appelle \og voisinage tubulaire \fg{} de $Y$ une famille $\Ucal=(U_w)_{w}$, où $w$ parcourt les places de $\Kb$ et $U_w$ est un voisinage ouvert de $Y (\overline{K}_w)$ dans la topologie $w$-adique, ces voisinages étant d'une certaine manière uniformes. L'exemple le plus naturel est la donnée d'un ouvert $U$ de $X(\C)$ contenant $Y(\C)$ (et alors les $U_w$ pour $w$ archimédienne seront tous égaux à $U$) et $U_w$ l'ensemble des points se réduisant dans $Y$ modulo $w$ (dans le modèle $\Xcal$) pour $w$ non archimédienne. On peut aussi voir un voisinage tubulaire comme une famille d'ouverts définis par une distance arithmétique à $Y$ bornée, au sens de Vojta (\cite{Vojtadiophapp}, paragraphe 2.5). On suppose précisément que quelle que soit la place $w$, notre point $P$ de $(\Xcal \backslash \Dcal) (\Ocal_{L,S_L})$ n'appartient jamais à $U_w$.  Alors, en reprenant le principe de la méthode de Runge, la condition de Runge devient 
\begin{equation}
\label{eqcondRungetub}
m_Y |S_L| < r,
\end{equation}
avec $m_Y$ le plus grand entier tel que l'intersection de $m_Y$ diviseurs parmi $D_1, \cdots, D_r$ n'est pas incluse dans $Y$. Si les zones de grande intersection des diviseurs sont contenues dans $Y$, le nombre $m_Y$ est donc bien plus petit que $m$. Ceci rend la condition \eqref{eqcondRungetub} satisfaisable et donne le \og théorème de Runge tubulaire \fg{} de l'article dont voici un énoncé simplifié (qui se généralise à $\Xcal$ normal, $D_1, \cdots, D_r$ de Cartier et définis sur une extension de $K$, et éventuellement gros au lieu d'amples).

\begin{thm}[\cite{LeFourn3}, Théorème 1 simplifié]
\hspace*{\fill}

	Avec les notations précédentes, si $D_1, \cdots, D_r$ sont amples, alors 
	
\[
	\bigcup_{\substack{(L,S_L) \\ m_Y |S_L|<r}} (\Xcal  \backslash \Dcal) (\Ocal_{L,S_L}) \mathlarger{\mathlarger{\backslash}} \bigcup_{w} U_w \textrm{ est fini,}
\]
et ce pour tout choix de voisinage tubulaire $\Ucal = (U_w)_{w}$ de $Y$.
\end{thm}

La preuve rigoureuse (qui suit ces idées) passe par une traduction de l'intégralité en termes de fonctions auxiliaires, puis un résultat de type Nullstellensatz pour ces fonctions. 

Notons que contrairement au cas des courbes, la finitude d'un ensemble de points entiers n'est pas garantie en général, c'est donc déj\`a un résultat non trivial à $(L,S_L)$ fixé. Ensuite, le fonctionnement de la preuve indique implicitement qu'elle peut encore une fois être adaptée en une méthode en pratique. Cela est un avantage sur d'autres résultats de finitude tels que par exemple le théorème CLZ de \cite{CorvajaLevinZannier} (basé sur le théorème du sous-espace de Schmidt), qui lui aussi fonctionne en excluant un fermé $Y$ (mais pour toute paire $(L,S_L)$).

Nous allons maintenant exposer l'application de Runge tubulaire aux variétés modulaires de Siegel $A_2(n)$.

\section{\texorpdfstring{Runge pour $A_2(n)$}{Runge pour A2(n)}}
\label{secvarmod}
Pour $n \geq 1$, la variété $A_2(n)$ sur $\Q(\zeta_n)$ est la variété modulaire paramétrant les classes d'isomorphismes de triplets $(A,\lambda,\alpha_n)$ où $A$ est une surface abélienne, $\lambda$ une polarisation principale de $A$ et $\alpha_n$ une structure symplectique de niveau $n$ pour $(A,\lambda)$ (autrement dit une base de la $n$-torsion de $A$ qui est symplectique pour l'accouplement de Weil). On note $A_2(n)^S$ sa compactification de Satake, qui est une variété normale projective de dimension 3 sur $\Q(\zeta_n)$, et $\partial A_2(n)^S:= A_2(n)^S \backslash A_2(n)$ son bord (de dimension 1). Il existe des modèles normaux naturels de $A_2(n)$ et $A_2(n)^S$ sur $\Z[\zeta_n,1/n]$ qu'on note respectivement $\Acal_2(n)$ et $\Acal_2(n)^S$ (\cite{ChaiFaltings}, chapitres IV et V).

On souhaite définir une bonne notion de points entiers sur $A_2(n)$, le problème étant qu'on ne dispose pas immédiatement de diviseurs agréables à choisir. Une première approche est d'exclure le bord, mais celui-ci est de codimension 2. Une seconde est de considérer les diviseurs paramétrant les produits de courbes elliptiques, mais à part pour $n=1$ ou $2$, ceux-ci ne sont pas amples. En fait, l'amplitude d'un diviseur de Cartier de $A_2(n)$ pour $n$ quelconque est dès $n \geq 2$ un problème difficile (on ne sait pas encore décrire convenablement le groupe de Picard de cette variété), mais dans un cas elle est immédiate : on va considérer les diviseurs des zéros des fonctions thêta sur $A_2(n)$, lorsque $n$ est pair. 

Rappelons qu'étant donné une variété abélienne principalement polarisée avec structure de niveau 2, on peut définir son diviseur thêta comme le diviseur des zéros de l'unique section non nulle du fibré ample symétrique induisant la polarisation (ce fibré est uniquement déterminé par le choix de 2-structure pour une normalisation dûe à Igusa que nous ne détaillerons pas). Les diviseurs que nous considérerons sur $A_2(n)$, notés les $D_{n,a,b}$ ($(a,b) \in (\Z/n\Z)^4$), seront alors les diviseurs paramétrant les surfaces abéliennes principalement polarisées $(A,\lambda)$ telles que le point de coordonnées $(a,b)$ dans la base symplectique $\alpha_n$ de la $n$-torsion appartient au diviseur thêta. En fait, ces diviseurs sont exactement les diviseurs de zéros de certaines fonctions thêta classique, d'où on déduit aisément leur amplitude. Par symétrie, $D_{n,a,b} = D_{n,-a,-b}$ et ce sont vraiment des diviseurs sauf pour 6 couples $(a,b)$ possibles qui correspondent aux six points de 2-torsion de $A$ se situant toujours sur le diviseur thêta (d'où $n^4/2 + 2$ diviseurs distincts). Les points de torsion se situant éventuellement sur le diviseur thêta et différents de ces six points seront dits non triviaux.

Le premier théorème sur les points entiers de $A_2(n)$ est alors le suivant.

\begin{thm}[\cite{LeFourn3}, Théorème 2]
Soit $n \geq 2$ pair.
Pour $U$ un voisinage ouvert de $\partial A_2(n)^S (\C)$ dans $A_2(n)^S(\C)$, on note $\Ecal(U)$ l'ensemble des points $P = \overline{(A,\lambda,\alpha_n)} \in A_2(n)(\Qb)$ tels que (si $L \supset \Q(\zeta_n)$ désigne un corps de définition de $P$)

\begin{itemize}
	\item La surface abélienne $A$ a potentiellement bonne réduction en toute place finie de $L$.
	\item Pour tout plongement $\sigma : L \rightarrow \C$, $P_\sigma \in A_2(n)(\C)$ n'appartient pas à $U$.
	\item Le nombre $s_L$ de places $v$ de $L$ telles que $v|\infty, v|2$, ou la réduction modulo $v$ de $P$ admet un point de $n$-torsion non trivial dans son diviseur thêta vérifie 
	\[
		(n^2 - 3) s_L < n^4/2 + 2.
	\]
Alors, quel que soit le choix de voisinage ouvert $U$, l'ensemble $\Ecal(U)$ est fini.
\end{itemize}
\end{thm}

Pour appréhender ce théorème, il est bon de penser au cas $n=2$, car alors les diviseurs $D_{2,a,b}$ se trouvent être exactement les dix diviseurs irréductibles de $A_2(2)$ paramétrant les produits de courbes elliptiques et la condition devient $s_L<10$. La potentielle bonne réduction signifie la réduction du point $P$ de $A_2(n)$ hors du bord ( c'est donc la condition tubulaire non-archimédienne), la seconde est la condition tubulaire archimédienne et les mauvaises places sont celles indiquées, auxquelles on ajoute les places au-dessus de $n$ car le modèle $\Acal_2(n)$ n'y est pas bien défini. Ces dernières sont au nombre de $(n/2) [L:\Q(\zeta_n)]$ au maximum, ce qui laisse de la marge pour les couples $(L,S_L)$.

La preuve de ce théorème est donc une application du théorème de Runge, le seul point restant étant l'évaluation de l'entier $m_{\partial A_2(n)}$, qui est égal à $(n^2 -3)$, borne d'ailleurs atteinte en les produits de courbes elliptiques. La condition tubulaire est absolument essentielle, l'entier $m$ (sans choix de fermé donc) étant asymptotiquement de l'ordre de $n^4/2$, soit le nombre total de diviseurs.

Nous allons maintenant passer au cas $n=2$, pour lequel la variété $A_2(2)$ et les fonctions thêta sont suffisamment bien connues pour avoir un résultat explicite de finitude.

\section{\texorpdfstring{Le cas $A_2(2)$}{Le cas A2(2)}}
\label{seccasdeux}

Pour une famille concrète de voisinages ouverts du bord dans $A_2(2)$, nous avons choisi de référer au modèle complexe de cette variété, à savoir le demi-espace supérieur de Siegel $\Hcal_2$ quotienté par le sous-groupe $\Gamma_2(2)$ du groupe symplectique entier de degré 4. Dans le domaine fondamental $\Fcal$ du quotient $\Hcal_2 / \Gamma_2(1)$ (cf. \cite{Klingen}, section I.2), pour un réel $t  \geq \sqrt{3}/2$, les éléments appartenant à l'ouvert $U_t$ sont les matrices $\tau = \begin{pmatrix} \tau_1 & \tau_2 \\ \tau_2 & \tau_4 \end{pmatrix}$ telles que $\Im (\tau_4) \geq t$. Par abus de notation, on renote $U_t \subset A_2(2)^S(\C)$ l'ensemble des points du bord ou tels que l'unique représentant de la surface abélienne paramétrée dans $\Fcal$ appartient à $U_t$ : c'est bien un voisinage ouvert de $\partial A_2(2)^S(\C)$.

Un autre point important est que les fonctions thêta définissent un plongement projectif $\psi$ dans $\P^9$. En effet, pour $m=(a,b) \in \{0,1/2\}^4$, définissons la fonction $\Theta_m$ sur $\Hcal_2$ par 
\[
	\Theta_m (\tau) = \sum_{n \in \Z^2} e^{2 i \pi {}^t (n+a) \tau (n+a) + 2 i \pi n {}^t b }.
\]
Les puissances quatrièmes de ces fonctions sont des formes modulaires de poids 2 pour $\Gamma_2(2)$. Parmi elles, six sont identiquement nulles (correspondant aux six choix de $(a,b)$ triviaux déj\`a mentionnés), on indexe par $E$ les dix autres. Alors, le morphisme 
\[
	\fonction{\psi}{A_2(2)(\C)}{\P^9(\C)}{\overline{\tau}}{(\Theta_m^4(\tau))_{m \in E}}
\]
est bien défini et induit un plongement projectif de $A_2(2)^S(\C)$ dans $\P^9$. De plus, il existe des relations linéaires entre les coordonnées qui permettent de voir  $\psi(A_2(2)^S(\C))$ comme une quartique de $\P^5$ ((\cite{vdG82}, Théorème 5.2) et par des résultats de Pazuki \cite{Pazuki12}, on peut lier la hauteur de $\psi(P)$ à la hauteur de Faltings stable de la surface $A$ que représente $P$, notée $h_\Fcal(A)$. On obtient alors le résultat suivant.

\begin{thm}[Version explicite pour $n=2$]
Soit $K$ un corps de nombres et $P=\overline{(A,\lambda,\alpha_2)} \in A_2(2)(K)$ où $A$ a potentiellement bonne réduction en toute place finie.

Soit $s_P$ le nombre d'idéaux premiers $\gP$ de $\Ocal_K$ tels que la réduction semistable de $A$ modulo un (tout) idéal au-dessus de $\gP$ est un produit de courbes elliptiques. Alors : 

$(a)$ Si $K=\Q$ ou un corps quadratique imaginaire et $s_P<4$, alors
\[
	h(\psi(P)) \leq 10.75, \quad h_\Fcal(A) \leq 1070.
\]

$(b)$ Soit $t \geq \sqrt{3}/2$. Si pour tout plongement $\sigma : K \rightarrow \C$, le point $P_\sigma \in A_2(2)_\C$ n'appartient pas à $U_t$ et
 \[
s_P + \left|\textrm{places archimédiennes de }K\right| < 10
\]
alors
\[
                                                                                                   h(\psi(P)) \leq 4 \pi t + 6.14, \quad h_\Fcal(A) \leq 2 \pi t + 535 \log(2 \pi t + 9).           		
\]
\end{thm}

La preuve de ce résultat effectif repose sur des estimations de $|\psi(P)|_v$ pour toute place $v$. Lorsque $v$ est archimédienne, on utilise les développements de Fourier des fonctions thêta pour montrer que six coordonnées tendent vers 0 lorsqu'on s'approche du bord et seulement une peut être proche de 0 loin du bord. Lorsque $v$ est nonarchimédienne et pas au-dessus de 2, on sait que les fonctions thêta s'algébrisent et passent à la réduction, en particulier si $A$ a potentiellement bonne réduction, la réduction modulo $v$ de $\psi(P)$ ne peut avoir qu'une seule coordonnée nulle et cela seulement quand $A$ modulo $v$ est un produit de courbes elliptiques. Enfin, lorsque $v$ divise 2, on ne peut plus utiliser cette algébrisation, mais on se sert des invariants d'Igusa pour les jacobiennes, qui décrivent eux aussi la forme de la réduction selon le théorème 1 de \cite{Liu93}. Il faut lier ces invariants à des formes modulaires puis aux fonctions thêta elles-mêmes, ce qui a exigé des calculs sur Sage \footnote{voir \\ \url{http://perso.ens-lyon.fr/samuel.le_fourn/contenu/fichiers_publis/Igusainvariants.ipynb}} pour assurer des estimations finales explicites. La raison de l'absence d'hypothèse archimédienne pour le $(a)$ est que dans cette situation, on peut ne rien supposer pour l'unique place archimédienne, ce qui élimine six coordonnées de $\psi(P)$, et il doit en rester une vu le principe de la méthode d'où la condition $s_P<4$.

Pour éclairer le lecteur, cette situation où la condition tubulaire archimédienne peut être enlevée ne marchera que pour $A_2(2)$, mais un résultat de la forme de $(b)$ est en théorie parfaitement envisageable pour tout $n$ pair, avec des calculs similaires. Il faudra cependant surmonter certaines difficultés, comme l'identification plus précise des zéros des fonctions thêta en jeu, qui du point de vue analytique revient à des estimations de développements de Fourier et du point de vue modulaire consiste en une meilleure compréhension des \og causes structurelles \fg{} de la présence d'un point de $n$-torsion non trivial dans le diviseur thêta d'une surface abélienne. Des applications de ces résultats à des surfaces modulaires de Hilbert ou des courbes de Shimura ne sont pas non plus à exclure, pour celles qui se plongent dans un $A_2(n)$.

\bibliographystyle{alphaSLF}
\bibliography{bibliotdn}
\end{document}